\def\marker{\>\hbox{${\vcenter{\vbox{
    \hrule height 0.4pt\hbox{\vrule width 0.4pt height 6pt
    \kern6pt\vrule width 0.4pt}\hrule height 0.4pt}}}$}\>}
\def\gpic#1{#1
     \medskip\par\noindent{\centerline{\box\graph}} \medskip}
\documentclass[11pt]{article}
\usepackage{fullpage}
\usepackage{jeffe}
\newtheorem{theorem}{Theorem}
\begin{document}
\author{Daniel W. Cranston\footnote{University of Illinois, Urbana-Champaign} \\ \texttt{dcransto@uiuc.edu}}
\title{Nomadic Decompositions of Bidirected Complete Graphs}
\maketitle

\abstract{
We use $K^*_n$ to denote the bidirected complete graph on $n$ vertices.
A nomadic Hamiltonian decomposition of $K^*_n$ is a Hamiltonian decomposition,
with the additional property that ``nomads'' walk along the Hamiltonian cycles (moving one vertex per time step) without colliding.
A nomadic near-Hamiltonian decomposition is defined similarly, except that the cycles in the decomposition have length
$n-1$, rather than length $n$.  J.A. Bondy asked whether these decompositions of $K^*_n$ exist for all $n$.  
We show that $K^*_n$ admits a nomadic near-Hamiltonian decomposition when $n\not\equiv 2\bmod 4$.
}
\bigskip

\textbf{Keywords:} Hamiltonian, decomposition, nomad
\section{Introduction}

A \textit{bidirected complete graph} on $n$ vertices, denoted $K^*_n$, is a digraph where each ordered 
pair of distinct vertices forms an edge.
Bermond and Faber first posed the following question: Can we partition the edge set of $K^*_n$ into Hamiltonian cycles?
Kirkman~\cite{berge} knew that this is possible when the number of vertices, $n$, is odd.
For even $n$, Tillson~\cite{tillson} showed this is possible for $n\geq 8$ (it cannot be done for smaller $n$).

At the 4th Cracow Conference at Czorstyn, Poland, (September 16-20, 2002)
J.A. Bondy~\cite{bondy} asked the following stronger version of the question.
Can we partition the edge set of $K^*_n$ into Hamiltonian cycles and put a ``nomad'' on each cycle so that if each
nomad moves forward one vertex (along his cycle) at each time step, then no two nomads ever collide?  
More formally, let $C_1, C_2, \ldots, C_{n-1}$ be a parition of $K^*_n$ into Hamiltonian cycles.
For $v\in V(C_i)$, let $v^{+k}$ denote the vertex we reach by starting at $v$ and following $C_i$ for $k$ steps.
Let $f$ be a function that chooses for each $C_i$ a root vertex $v_i$.
If there exists a function $f$ such that $v_i^{+k}\neq v_j^{+k}$ for all $i\neq j$ and for all $k$,
then we say that the decomposition of $K^*_n$ is a \textit{nomadic Hamiltonian decomposition}.

We call a (directed) cycle of length $n-1$ a near-Hamiltonian cycle.
Analagous to a nomadic Hamiltonian decomposition, we define a \textit{nomadic near-Hamiltonian decomposition}.
Bondy also asked if $K^*_n$ has a nomadic near-Hamiltonian decomposition for every $n$.
In this paper, we show that $K^*_n$ does have a nomadic near-Hamiltonian decomposition for every $n\not\equiv 2\bmod 4$.

\section{Nomadic Near-Hamiltonian Decomposition of $K^*_n$ for $n$ odd}

\begin{figure}[!thb]
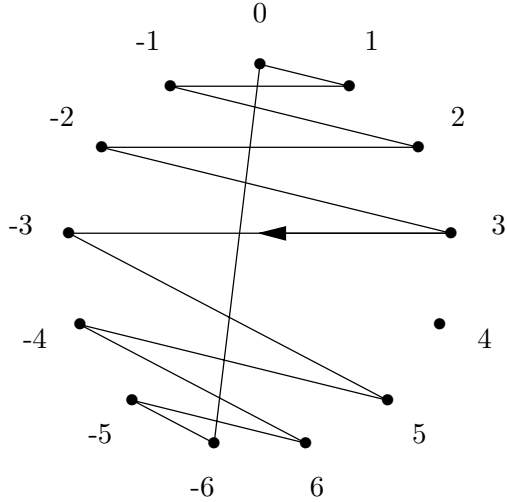

\gpic{
\expandafter\ifx\csname graph\endcsname\relax \csname newbox\endcsname\graph\fi
\expandafter\ifx\csname graphtemp\endcsname\relax \csname newdimen\endcsname\graphtemp\fi
\setbox\graph=\vtop{\vskip 0pt\hbox{%
    \graphtemp=.5ex\advance\graphtemp by 1.616in
    \rlap{\kern 0.308in\lower\graphtemp\hbox to 0pt{\hss $\bullet$\hss}}%
    \graphtemp=.5ex\advance\graphtemp by 1.138in
    \rlap{\kern 0.250in\lower\graphtemp\hbox to 0pt{\hss $\bullet$\hss}}%
    \graphtemp=.5ex\advance\graphtemp by 0.687in
    \rlap{\kern 0.421in\lower\graphtemp\hbox to 0pt{\hss $\bullet$\hss}}%
    \graphtemp=.5ex\advance\graphtemp by 0.368in
    \rlap{\kern 0.782in\lower\graphtemp\hbox to 0pt{\hss $\bullet$\hss}}%
    \graphtemp=.5ex\advance\graphtemp by 0.252in
    \rlap{\kern 1.250in\lower\graphtemp\hbox to 0pt{\hss $\bullet$\hss}}%
    \graphtemp=.5ex\advance\graphtemp by 0.368in
    \rlap{\kern 1.718in\lower\graphtemp\hbox to 0pt{\hss $\bullet$\hss}}%
    \graphtemp=.5ex\advance\graphtemp by 0.687in
    \rlap{\kern 2.079in\lower\graphtemp\hbox to 0pt{\hss $\bullet$\hss}}%
    \graphtemp=.5ex\advance\graphtemp by 1.138in
    \rlap{\kern 2.250in\lower\graphtemp\hbox to 0pt{\hss $\bullet$\hss}}%
    \graphtemp=.5ex\advance\graphtemp by 1.616in
    \rlap{\kern 2.192in\lower\graphtemp\hbox to 0pt{\hss $\bullet$\hss}}%
    \graphtemp=.5ex\advance\graphtemp by 2.013in
    \rlap{\kern 1.918in\lower\graphtemp\hbox to 0pt{\hss $\bullet$\hss}}%
    \graphtemp=.5ex\advance\graphtemp by 2.237in
    \rlap{\kern 1.490in\lower\graphtemp\hbox to 0pt{\hss $\bullet$\hss}}%
    \graphtemp=.5ex\advance\graphtemp by 2.237in
    \rlap{\kern 1.010in\lower\graphtemp\hbox to 0pt{\hss $\bullet$\hss}}%
    \graphtemp=.5ex\advance\graphtemp by 2.013in
    \rlap{\kern 0.582in\lower\graphtemp\hbox to 0pt{\hss $\bullet$\hss}}%
    \special{pn 8}%
    \special{pa 1250 252}%
    \special{pa 1718 368}%
    \special{pa 782 368}%
    \special{pa 2079 687}%
    \special{pa 421 687}%
    \special{pa 2250 1138}%
    \special{pa 250 1138}%
    \special{fp}%
    \special{pa 250 1138}%
    \special{pa 1918 2013}%
    \special{pa 308 1616}%
    \special{pa 1490 2237}%
    \special{pa 582 2013}%
    \special{pa 1010 2237}%
    \special{pa 1250 252}%
    \special{fp}%
    \graphtemp=.5ex\advance\graphtemp by 0.000in
    \rlap{\kern 1.250in\lower\graphtemp\hbox to 0pt{\hss 0\hss}}%
    \graphtemp=.5ex\advance\graphtemp by 0.145in
    \rlap{\kern 1.835in\lower\graphtemp\hbox to 0pt{\hss 1\hss}}%
    \graphtemp=.5ex\advance\graphtemp by 0.544in
    \rlap{\kern 2.286in\lower\graphtemp\hbox to 0pt{\hss 2\hss}}%
    \graphtemp=.5ex\advance\graphtemp by 1.108in
    \rlap{\kern 2.500in\lower\graphtemp\hbox to 0pt{\hss 3\hss}}%
    \graphtemp=.5ex\advance\graphtemp by 1.706in
    \rlap{\kern 2.427in\lower\graphtemp\hbox to 0pt{\hss 4\hss}}%
    \graphtemp=.5ex\advance\graphtemp by 2.202in
    \rlap{\kern 2.085in\lower\graphtemp\hbox to 0pt{\hss 5\hss}}%
    \graphtemp=.5ex\advance\graphtemp by 2.481in
    \rlap{\kern 1.550in\lower\graphtemp\hbox to 0pt{\hss 6\hss}}%
    \graphtemp=.5ex\advance\graphtemp by 2.481in
    \rlap{\kern 0.950in\lower\graphtemp\hbox to 0pt{\hss -6\hss}}%
    \graphtemp=.5ex\advance\graphtemp by 2.202in
    \rlap{\kern 0.415in\lower\graphtemp\hbox to 0pt{\hss -5\hss}}%
    \graphtemp=.5ex\advance\graphtemp by 1.706in
    \rlap{\kern 0.073in\lower\graphtemp\hbox to 0pt{\hss -4\hss}}%
    \graphtemp=.5ex\advance\graphtemp by 1.108in
    \rlap{\kern 0.000in\lower\graphtemp\hbox to 0pt{\hss -3\hss}}%
    \graphtemp=.5ex\advance\graphtemp by 0.544in
    \rlap{\kern 0.214in\lower\graphtemp\hbox to 0pt{\hss -2\hss}}%
    \graphtemp=.5ex\advance\graphtemp by 0.145in
    \rlap{\kern 0.665in\lower\graphtemp\hbox to 0pt{\hss -1\hss}}%
    \special{pa 2250 1138}%
    \special{pa 1250 1138}%
    \special{fp}%
    \special{sh 1.000}%
    \special{pa 1384 1172}%
    \special{pa 1250 1138}%
    \special{pa 1384 1104}%
    \special{pa 1384 1172}%
    \special{fp}%
    \hbox{\vrule depth2.481in width0pt height 0pt}%
    \kern 2.500in
  }%
}%
}
\caption{A single near-Hamiltonian cycle in which each edge has distinct length.  
A nomad who begins at vertex 0 follows edges with successive lengths $1,-2,3,-4,5,-6,-5,4,-3,2,-1,6$.
We can partition the edges of $K^*_{13}$ into $13$ copies of this cycle.
Together, these 13 cycles form an r.s. nomadic near-Hamiltonian decomposition of $K^*_{13}$.
}
\end{figure}

Label the vertices of $K^*_n$ with the integers $0, 1, \ldots, n-1$.
Let the positions of two nomads be given by $g(t)$ and $h(t)$.
We say the nomads are \textit{rotationally symmetric} (or r.s.) if there exists a constant $c\neq 0$ 
such that for every $t$ we have $g(t)-h(t) = c$.
Rotationally symmetric nomads are convenient to study because they never occupy the same vertex.  
Since rotational symmetry is transitive, we can ask the stronger question:
Does $K^*_n$ have a nomadic near-Hamiltonian decomposition that consists entirely of nomads that are rotationally symmetric?
We call this an r.s. nomadic near-Hamiltonian decomposition.

Let $(g(t) - g(t-1)) \bmod n$ be the \textit{length} of the edge that nomad $g$ walks along at time $t$.
Since there are no edges of length $0$, there are exactly $n-1$ different possible edge lengths.
In $K^*_n$ there are $n$ edges of each length, so if $n$ cycles are rotationally symmetric, 
then each cycle must contain exactly one edge of each length.  
Hence, to find an r.s. nomadic near-Hamiltonian decomposition, 
it is sufficient to find a single near-Hamiltonian cycle in which each edge has distinct length.
This leads to the following theorem.

\begin{theorem}
If $n$ is odd, then $K^*_n$ has an r.s. nomadic near-Hamiltonian decomposition.
\end{theorem}
\begin{proof}
To prove the theorem, we only need to find a single near-Hamiltonian cycle in which each edge has distinct length.
If we find such a cycle, we can place $n$ rotationally symmetric nomads, one of whom walks along this cycle.
Since $n$ is odd, let $k = (n-1)/2$.
It is convenient to 
label the vertices of the cycle with the integers $-k, -(k-1), \ldots, -1, 0, 1, \ldots, k-1, k$.
One nomad starts at vertex 0 at time 0.  
At time $t$, his position is the sum (modulo $n$) of the first $t$ terms of the sequence 
$1, -2, 3, -4, \ldots \pm k, \pm (k-1), \mp (k-2), \ldots, -3, 2, -1, \mp k$.

Note that the length of the edge the nomad follows at time $t$ is the $t$th term in the above sequence.
Since the terms of the sequence are distinct, so are the edge lengths of the cycle.  Hence, we only need to verify that
the nomad follows a cycle of length $n-1$ (i.e. he doesn't revisit any vertex too soon).

Two time steps after the nomad is at vertex $1\leq i < k$, he is at vertex $i+1$;
similarly, two time steps after the nomad is at vertex $-(k-1) < i \leq 0$, he is at vertex $i-1$.
The single exception is when the nomad follows edges of length $\pm k$ and $\pm (k-1)$ in succession.
At this point, the nomad ``skips over'' vertex $(-1)^{\frac{n-1}2}\ceil{\frac{n+1}4}$.
Verifying that all is well when the nomad leaves vertices $-(k-1), -k$, and $k$ requires only a short case analysis.

\end{proof}

\section{Nomadic Near-Hamiltonian Decomposition of $K^*_n$ for $n\equiv 0\bmod 4$}

It is natural to try the idea we used for odd $n$ to construct a nomadic near-Hamiltonian decomposition for even $n$.
However, this approach is doomed to fail.  
In an r.s. near-Hamiltonian decomposition, each cycle must contain one edge of each length.  
Thus, the ``length'' of the cycle must be $\sum_{i=1}^{n-1}i = n(n-1)/2$. However $ n(n-1)/2 \not\equiv 0 \bmod n$.  
In other words, when $n$ is even, a nomad who follows a path consisting of one edge of each length 
won't end up back where he started.
In spite of this, when $n\equiv 0\bmod 4$, we are able to construct a nomadic near-Hamiltonian 
decomposition using two disjoint sets of rotationally symmetric nomads (we call these sets $A$ and $B$).

\begin{figure}[!thb]
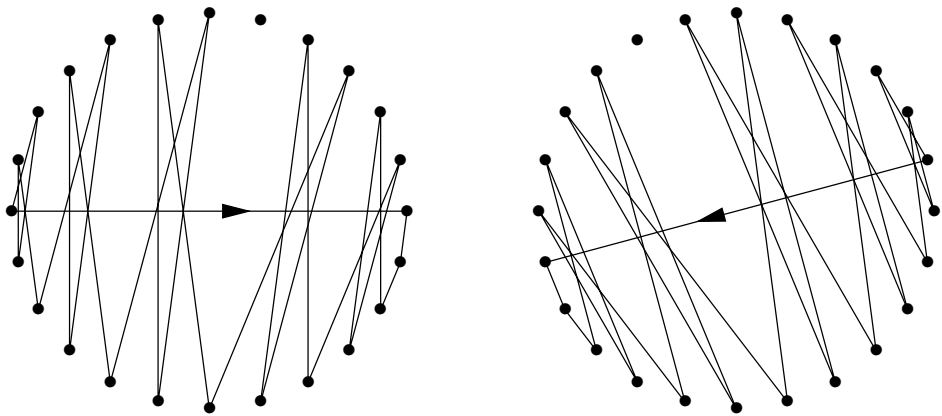

\gpic{
\expandafter\ifx\csname graph\endcsname\relax \csname newbox\endcsname\graph\fi
\expandafter\ifx\csname graphtemp\endcsname\relax \csname newdimen\endcsname\graphtemp\fi
\setbox\graph=\vtop{\vskip 0pt\hbox{%
    \graphtemp=.5ex\advance\graphtemp by 1.852in
    \rlap{\kern 1.852in\lower\graphtemp\hbox to 0pt{\hss $\bullet$\hss}}%
    \graphtemp=.5ex\advance\graphtemp by 2.017in
    \rlap{\kern 1.638in\lower\graphtemp\hbox to 0pt{\hss $\bullet$\hss}}%
    \graphtemp=.5ex\advance\graphtemp by 2.120in
    \rlap{\kern 1.389in\lower\graphtemp\hbox to 0pt{\hss $\bullet$\hss}}%
    \graphtemp=.5ex\advance\graphtemp by 2.155in
    \rlap{\kern 1.121in\lower\graphtemp\hbox to 0pt{\hss $\bullet$\hss}}%
    \graphtemp=.5ex\advance\graphtemp by 2.120in
    \rlap{\kern 0.853in\lower\graphtemp\hbox to 0pt{\hss $\bullet$\hss}}%
    \graphtemp=.5ex\advance\graphtemp by 2.017in
    \rlap{\kern 0.603in\lower\graphtemp\hbox to 0pt{\hss $\bullet$\hss}}%
    \graphtemp=.5ex\advance\graphtemp by 1.852in
    \rlap{\kern 0.389in\lower\graphtemp\hbox to 0pt{\hss $\bullet$\hss}}%
    \graphtemp=.5ex\advance\graphtemp by 1.638in
    \rlap{\kern 0.225in\lower\graphtemp\hbox to 0pt{\hss $\bullet$\hss}}%
    \graphtemp=.5ex\advance\graphtemp by 1.389in
    \rlap{\kern 0.121in\lower\graphtemp\hbox to 0pt{\hss $\bullet$\hss}}%
    \graphtemp=.5ex\advance\graphtemp by 1.121in
    \rlap{\kern 0.086in\lower\graphtemp\hbox to 0pt{\hss $\bullet$\hss}}%
    \graphtemp=.5ex\advance\graphtemp by 0.853in
    \rlap{\kern 0.121in\lower\graphtemp\hbox to 0pt{\hss $\bullet$\hss}}%
    \graphtemp=.5ex\advance\graphtemp by 0.603in
    \rlap{\kern 0.225in\lower\graphtemp\hbox to 0pt{\hss $\bullet$\hss}}%
    \graphtemp=.5ex\advance\graphtemp by 0.389in
    \rlap{\kern 0.389in\lower\graphtemp\hbox to 0pt{\hss $\bullet$\hss}}%
    \graphtemp=.5ex\advance\graphtemp by 0.225in
    \rlap{\kern 0.603in\lower\graphtemp\hbox to 0pt{\hss $\bullet$\hss}}%
    \graphtemp=.5ex\advance\graphtemp by 0.121in
    \rlap{\kern 0.853in\lower\graphtemp\hbox to 0pt{\hss $\bullet$\hss}}%
    \graphtemp=.5ex\advance\graphtemp by 0.086in
    \rlap{\kern 1.121in\lower\graphtemp\hbox to 0pt{\hss $\bullet$\hss}}%
    \graphtemp=.5ex\advance\graphtemp by 0.121in
    \rlap{\kern 1.389in\lower\graphtemp\hbox to 0pt{\hss $\bullet$\hss}}%
    \graphtemp=.5ex\advance\graphtemp by 0.225in
    \rlap{\kern 1.638in\lower\graphtemp\hbox to 0pt{\hss $\bullet$\hss}}%
    \graphtemp=.5ex\advance\graphtemp by 0.389in
    \rlap{\kern 1.852in\lower\graphtemp\hbox to 0pt{\hss $\bullet$\hss}}%
    \graphtemp=.5ex\advance\graphtemp by 0.603in
    \rlap{\kern 2.017in\lower\graphtemp\hbox to 0pt{\hss $\bullet$\hss}}%
    \graphtemp=.5ex\advance\graphtemp by 0.853in
    \rlap{\kern 2.120in\lower\graphtemp\hbox to 0pt{\hss $\bullet$\hss}}%
    \graphtemp=.5ex\advance\graphtemp by 1.121in
    \rlap{\kern 2.155in\lower\graphtemp\hbox to 0pt{\hss $\bullet$\hss}}%
    \graphtemp=.5ex\advance\graphtemp by 1.389in
    \rlap{\kern 2.120in\lower\graphtemp\hbox to 0pt{\hss $\bullet$\hss}}%
    \graphtemp=.5ex\advance\graphtemp by 1.638in
    \rlap{\kern 2.017in\lower\graphtemp\hbox to 0pt{\hss $\bullet$\hss}}%
    \special{pn 8}%
    \special{pa 2155 1121}%
    \special{pa 2120 1389}%
    \special{pa 2017 1638}%
    \special{pa 2017 603}%
    \special{pa 1852 1852}%
    \special{pa 2120 853}%
    \special{pa 1638 2017}%
    \special{pa 1638 225}%
    \special{pa 1389 2120}%
    \special{fp}%
    \special{pa 1389 2120}%
    \special{pa 1852 389}%
    \special{pa 1121 2155}%
    \special{pa 853 121}%
    \special{pa 853 2120}%
    \special{pa 1121 86}%
    \special{pa 603 2017}%
    \special{pa 389 389}%
    \special{pa 389 1852}%
    \special{fp}%
    \special{pa 389 1852}%
    \special{pa 603 225}%
    \special{pa 225 1638}%
    \special{pa 121 853}%
    \special{pa 121 1389}%
    \special{pa 225 603}%
    \special{pa 86 1121}%
    \special{fp}%
    \special{pa 86 1121}%
    \special{pa 1328 1121}%
    \special{fp}%
    \special{sh 1.000}%
    \special{pa 1190 1086}%
    \special{pa 1328 1121}%
    \special{pa 1190 1155}%
    \special{pa 1190 1086}%
    \special{fp}%
    \special{pa 1328 1121}%
    \special{pa 2155 1121}%
    \special{fp}%
    \graphtemp=.5ex\advance\graphtemp by 1.852in
    \rlap{\kern 4.611in\lower\graphtemp\hbox to 0pt{\hss $\bullet$\hss}}%
    \graphtemp=.5ex\advance\graphtemp by 2.017in
    \rlap{\kern 4.397in\lower\graphtemp\hbox to 0pt{\hss $\bullet$\hss}}%
    \graphtemp=.5ex\advance\graphtemp by 2.120in
    \rlap{\kern 4.147in\lower\graphtemp\hbox to 0pt{\hss $\bullet$\hss}}%
    \graphtemp=.5ex\advance\graphtemp by 2.155in
    \rlap{\kern 3.879in\lower\graphtemp\hbox to 0pt{\hss $\bullet$\hss}}%
    \graphtemp=.5ex\advance\graphtemp by 2.120in
    \rlap{\kern 3.611in\lower\graphtemp\hbox to 0pt{\hss $\bullet$\hss}}%
    \graphtemp=.5ex\advance\graphtemp by 2.017in
    \rlap{\kern 3.362in\lower\graphtemp\hbox to 0pt{\hss $\bullet$\hss}}%
    \graphtemp=.5ex\advance\graphtemp by 1.852in
    \rlap{\kern 3.148in\lower\graphtemp\hbox to 0pt{\hss $\bullet$\hss}}%
    \graphtemp=.5ex\advance\graphtemp by 1.638in
    \rlap{\kern 2.983in\lower\graphtemp\hbox to 0pt{\hss $\bullet$\hss}}%
    \graphtemp=.5ex\advance\graphtemp by 1.389in
    \rlap{\kern 2.880in\lower\graphtemp\hbox to 0pt{\hss $\bullet$\hss}}%
    \graphtemp=.5ex\advance\graphtemp by 1.121in
    \rlap{\kern 2.845in\lower\graphtemp\hbox to 0pt{\hss $\bullet$\hss}}%
    \graphtemp=.5ex\advance\graphtemp by 0.853in
    \rlap{\kern 2.880in\lower\graphtemp\hbox to 0pt{\hss $\bullet$\hss}}%
    \graphtemp=.5ex\advance\graphtemp by 0.603in
    \rlap{\kern 2.983in\lower\graphtemp\hbox to 0pt{\hss $\bullet$\hss}}%
    \graphtemp=.5ex\advance\graphtemp by 0.389in
    \rlap{\kern 3.148in\lower\graphtemp\hbox to 0pt{\hss $\bullet$\hss}}%
    \graphtemp=.5ex\advance\graphtemp by 0.225in
    \rlap{\kern 3.362in\lower\graphtemp\hbox to 0pt{\hss $\bullet$\hss}}%
    \graphtemp=.5ex\advance\graphtemp by 0.121in
    \rlap{\kern 3.611in\lower\graphtemp\hbox to 0pt{\hss $\bullet$\hss}}%
    \graphtemp=.5ex\advance\graphtemp by 0.086in
    \rlap{\kern 3.879in\lower\graphtemp\hbox to 0pt{\hss $\bullet$\hss}}%
    \graphtemp=.5ex\advance\graphtemp by 0.121in
    \rlap{\kern 4.147in\lower\graphtemp\hbox to 0pt{\hss $\bullet$\hss}}%
    \graphtemp=.5ex\advance\graphtemp by 0.225in
    \rlap{\kern 4.397in\lower\graphtemp\hbox to 0pt{\hss $\bullet$\hss}}%
    \graphtemp=.5ex\advance\graphtemp by 0.389in
    \rlap{\kern 4.611in\lower\graphtemp\hbox to 0pt{\hss $\bullet$\hss}}%
    \graphtemp=.5ex\advance\graphtemp by 0.603in
    \rlap{\kern 4.775in\lower\graphtemp\hbox to 0pt{\hss $\bullet$\hss}}%
    \graphtemp=.5ex\advance\graphtemp by 0.853in
    \rlap{\kern 4.879in\lower\graphtemp\hbox to 0pt{\hss $\bullet$\hss}}%
    \graphtemp=.5ex\advance\graphtemp by 1.121in
    \rlap{\kern 4.914in\lower\graphtemp\hbox to 0pt{\hss $\bullet$\hss}}%
    \graphtemp=.5ex\advance\graphtemp by 1.389in
    \rlap{\kern 4.879in\lower\graphtemp\hbox to 0pt{\hss $\bullet$\hss}}%
    \graphtemp=.5ex\advance\graphtemp by 1.638in
    \rlap{\kern 4.775in\lower\graphtemp\hbox to 0pt{\hss $\bullet$\hss}}%
    \special{pa 2880 1389}%
    \special{pa 2983 1638}%
    \special{pa 3148 1852}%
    \special{pa 2880 853}%
    \special{pa 3362 2017}%
    \special{pa 2845 1121}%
    \special{pa 3611 2120}%
    \special{pa 3148 389}%
    \special{pa 3879 2155}%
    \special{fp}%
    \special{pa 3879 2155}%
    \special{pa 2983 603}%
    \special{pa 4147 2120}%
    \special{pa 3879 86}%
    \special{pa 4397 2017}%
    \special{pa 3611 121}%
    \special{pa 4611 1852}%
    \special{pa 4397 225}%
    \special{pa 4775 1638}%
    \special{fp}%
    \special{pa 4775 1638}%
    \special{pa 4147 121}%
    \special{pa 4879 1389}%
    \special{pa 4775 603}%
    \special{pa 4914 1121}%
    \special{pa 4611 389}%
    \special{pa 4879 853}%
    \special{fp}%
    \special{pa 4879 853}%
    \special{pa 3679 1174}%
    \special{fp}%
    \special{sh 1.000}%
    \special{pa 3822 1172}%
    \special{pa 3679 1174}%
    \special{pa 3804 1105}%
    \special{pa 3822 1172}%
    \special{fp}%
    \special{pa 3679 1174}%
    \special{pa 2880 1389}%
    \special{fp}%
    \hbox{\vrule depth2.241in width0pt height 0pt}%
    \kern 5.000in
  }%
}%
}
\caption{
The cycle on the left is the route followed by a nomad in $A$.  
Each of the nomads in $A$ follows a route identical to this cycle, except that each other route is rotated by a distinct multiple of 2 vertices.
The cycle on the right is the route followed by a nomad in $B$.  
Each of the nomads in $B$ follows a route identical to this cycle, except that each other route is rotated by a distinct multiple of 2 vertices.
}
\end{figure}

\begin{theorem}
If $n\equiv 0\bmod 4$, then $K^*_n$ has a nomadic near-Hamiltonian decomposition.
\end{theorem}
\begin{proof}
Since $n$ is even, let $k=n/2$.
It is convenient to label the vertices with the integers $-(k-1), -(k-2), \ldots, -2, -1, 0, 1, 2, \ldots, {k-1}, k$.
Our $n$ nomads consist of two groups of $n/2$ rotationally symmetric nomads.
We call these groups $A$ and $B$.  
Group $A$ consists of the nomads that are initially on vertices with odd labels and
group $B$ consists of the nomads that are initially on vertices with even labels.
By definition, we won't have conflicts within $A$ or within $B$.
To avoid conflicts between $A$ and $B$,
we construct the cycles so that at each time step, either all nomads in $A$ are on vertices with odd labels and all nomads in 
$B$ are on vertices with even labels, or vice versa.  
To achieve this, it is sufficient if at each time step the nomads all follow edges with lengths of the same parity.

It is convenient to refer to the lengths of the edges as $-(k-1), -(k-2), \ldots, -2,-1, 1, 2,\ldots, k-1, k$.
Since the nomads within $A$ and within $B$ are rotationally symmetric, we only describe a single cycle in each of $A$ and $B$.
In fact, we construct the cycles in $A$ and $B$ so that at each time step the length of the edge in $B$ is 
negative the length of the edge in $A$.

The construction is best understood by an example.
Consider the cycles shown in Figure~2, where $k=12$.
The cycle on the left (which is in $A$) has successive edges of lengths 12, 1, 1, -4, 5, -4, 5, -8, 9, -8, 9, 11, -10, 11, -10, 7, -6, 7, -6, 3, -2, 3, -2.
The length of an edge at time $t$ in $B$ is negative the length of an edge at time $t$ in $A$.
Thus, the cycle on the right (which is in $B$) has successive edges of lengths 12, -1, -1, 4, -5, 4, -5, 8, -9, 8, -9, -11, 10, -11, 10, -7, 6, -7, 6, -3, 2, -3, 2.

To verify that the decomposition is valid, we must confirm two facts: 1) each edge of $K^*_n$ appears in exactly one cycle and 
2) each nomad follows a cycle of length $n-1$ (i.e. he doesn't revisit any vertex too soon).
The key idea in establishing the first fact is that if $A$ contains an edge of length $l$, 
then $A$ contains all the edges of length $l$, and $B$ contains all the edges of length $-l$; 
the single exception is that each cycle contains exactly one edge of length $k$.
Hence, for each $1\leq l\ < k$, each cycle either contains two edges of length $l$ and no edges of length $-l$, or it
contains two edges of length $-l$ and no edges of length $l$.
By examining the figures, it is easy to see that in the example each nomad follows a cycle of length $n-1$.
Since the length of an edge at time $t$ in $B$ is negative the length of an edge at time $t$ in $A$,
we only describe the general construction of a cycle in $A$.

To construct a cycle in $A$, we group successive edge lengths into blocks of four.
An \textit{increasing} block is of the form $-4i, 4i+1, -4i, 4i+1$
and a \textit{decreasing} block is of the form $4j-1, -4j+2, 4j-1, -4j+2$.
The list of edge lengths for a cycle in $A$ begins $k,1,1$, then is followed by $\floor{(n-4)/8}$ increasing blocks, 
and then by $\ceil{(n-4)/8}$ decreasing blocks.  
Among the increasing blocks, $i$ increases successively from 1 to $\floor{(n-4)/8}$, 
and among the decreasing blocks, $j$ decreases successively from $\ceil{(n-4)/8}$ to 1. 
In the example above, $i$ increases from 1 to 2, then $j$ decreases from 3 to 1.

In the general case, it is easiest to view a nomad in $A$ as steadily moving from right to left, alternating between vertices on
the ``top'' and ``bottom'' as he goes.
After following the four edges in a block (either increasing or decreasing), a nomad has moved two vertices clockwise.
Thus, the initial edge lengths $k,1,1$ move the nomad from left to right, then the increasing and decreasing blocks gradually
move the nomad back from right to left.
These observations make it easy to verify that each nomad follows a cycle of length $n-1$.
\end{proof}

\section{Discussion}
Although, we have not proved anything about the existence of nomadic Hamiltonian decompositions, we make one remark
concerning a stronger version of the question posed by Bondy.  
He asked if it is true for sufficiently large $n$ that \textit{every} Hamiltonian decomposition of $K^*_n$ is nomadic.  
We show the answer to this question is ``no.''  Let $n$ be prime.  We can
decompose $K^*_n$ into $n-1$ directed cycles, such that all of the edges within each cycle have the same length.  We show
that for this decomposition any two nomads will collide, regardless of their initial positions.  Say nomad 1 starts at
vertex $v_1$ on a cycle with edges of length $l_1$.  Similarly, say nomad 2 starts at vertex $v_2$ on a cycle with edges
of length $l_2$.  The two nomads will collide if and only if there exists time $t$ such that $t(l_2-l_1)\equiv v_1-v_2\bmod n$.
Since $n$ is prime, this equivalence does have a solution.

The problem of finding a nomadic near-Hamiltonian decomposition of $K^*_n$ remains open when $n\equiv 2\bmod 4$.
In fact, we have been unable to find such a decomposition for $n=6.$

We thank Douglas West for presenting this problem at the REGS combinatorics problem session at the University of Illinois, during summer 2006.
Thanks also to Tracy Grauman for many suggestions that improved the exposition.

\end{document}